\documentclass[12pt]{amsart}

\usepackage{amsmath,amsthm,amsfonts,amssymb,mathrsfs}

\newtheorem{theorem_etc}{theorem_etc}[section]
\newtheorem{theorem}[theorem_etc]{Theorem}

\newtheorem{lemma}[theorem_etc]{Lemma}
\newtheorem{corollary}[theorem_etc]{Corollary}
\newtheorem{definition}[theorem_etc]{Definition}
\newtheorem{remark}[theorem_etc]{Remark}

\newcommand{\N}{\mathbb{N}}

\newcommand{\R}{\mathbb{R}}
\newcommand{\Z}{\mathbb{Z}}

\newcommand{\cl}{\operatorname{\mathrm{cl}}}
\newcommand{\FP}{\operatorname{\mathit{FP}}}
\newcommand{\AP}{\operatorname{\mathit{AP}}}

\allowdisplaybreaks

\begin{document}

\title{On the topology of free paratopological groups}

\author{Ali Sayed Elfard}
\address{School of Mathematics and Applied Statistics,
University of Wollongong,
Wollongong, Australia}
\email{ae351@uow.edu.au\textrm{,} a.elfard@yahoo.com}

\author{Peter Nickolas}
\address{School of Mathematics and Applied Statistics,
University of Wollongong,
Wollongong, Australia}
\email{peter@uow.edu.au}

\keywords{Paratopological group, free paratopological group, quasi-pseudometric, quasi-prenorm, Graev extension, Joiner's lemma}

\subjclass[2000]{Primary 22A30; secondary 54D10, 54E99, 54H99}

\begin{abstract}
The result often known as Joiner's lemma is fundamental in understanding the topology of the free topological group $F(X)$ on a Tychonoff space~$X$. In this paper, an analogue of Joiner's lemma for the free paratopological group $\FP(X)$ on a $T_1$~space~$X$ is proved. Using this, it is shown that the following conditions are equivalent for a space~$X$: (1)~$X$ is~$T_1$; (2)~$\FP(X)$ is~$T_1$; (3)~the subspace $X$ of $\FP(X)$ is closed; (4)~the subspace $X^{-1}$ of $\FP(X)$ is discrete; (5)~the subspace $X^{-1}$ is~$T_1$; (6)~the subspace $X^{-1}$ is closed; and (7)~the subspace $\FP_n(X)$ is closed for all $n \in \N$, where $\FP_n(X)$ denotes the subspace of $\FP(X)$ consisting of all words of length at most~$n$.
\end{abstract}

\maketitle

\section{Introduction}

The notions of the free topological group on a Tychonoff space $X$ and a pointed Tychonoff space $(X, e)$ were introduced in the 1940s by Markov~\cite{Markov1,Markov2,Markov3} and Graev~\cite{Graev1,Graev2}, respectively. In both cases, the groups are Hausdorff. In 1976 Joiner~\cite{Joiner} provided a complete description of a neighbourhood basis at any word of length exactly~$n$ in the subspace $F_n(X)$ of the Graev free topological group on~$X$, where $F_n(X)$ denotes the set of all words in the group of length at most~$n$. Already in 1968 Arhangel'skii~\cite{Arh1} had proved essentially the same result as Joiner, though as noted in~\cite{Arh2} his result did not at the time attract much attention. Joiner's argument, though much more complex than that of Arhangel'skii (see~\cite{Arh2}), gives information not only about the topology of the free topological group but also about the topology induced on the free group by certain pseudometrics defined by Graev, and the result of Arhangel'skii and Joiner is commonly referred to as Joiner's lemma.

In 2003 Romaguera, Sanchis and Tkachenko \cite{RoSaTk} proved the existence of the free paratopological group $\FP(X, \mathcal{U})$ on a quasi-uniform space $(X, \mathcal{U})$ and investigated its separation properties. In 2006 Pyrch and Ravsky~\cite{PyRa} investigated some of the topological properties of the free paratopological group $\FP(X)$ on a topological space~$X$.

All of the authors above also discuss the corresponding free abelian topological or paratopological groups, and indeed some of the results of~\cite{PyRa} are proved in the abelian case only. For further background, the reader is referred to the introduction of~\cite{RoSaTk}, to Ravsky~\cite{Ravsky1,Ravsky2} and to Marin and Romaguera~\cite{MaRo}.

The main result of this paper, Theorem~\ref{joiner}, is an analogue of Joiner's lemma for free paratopological groups. The result takes the following form. Let $X$ be a $T_1$ space and denote by $\FP_n(X)$ the subspace of $\FP(X)$ consisting of all words of length at most~$n$. Suppose that $w = x_1^{\epsilon_1} x_2^{\epsilon_2} \ldots x_n^{\epsilon_n}$ is a reduced word in $\FP_n(X)$. Then a base at~$w$ in $\FP_n(X)$ is given by the collection of all sets of the form $U_1^{\epsilon_1} U_2^{\epsilon_2} \ldots U_n^{\epsilon_n}$, where for $i = 1, 2, \ldots, n$ the set $U_i$ is a neighbourhood of $x_i$ in $X$ when $\epsilon_i=1$ and $U_i = \{x_i\}$ when $\epsilon_i = -1$.

Using the above result and other ideas, we strengthen and generalise some results from Pyrch and Ravsky~\cite{PyRa} on the topological properties of $\FP(X)$ (see Theorems \ref{rez_thm}, \ref{Xpowers} and~\ref{equiv_conds}).

\section{Definitions and preliminaries}
\label{defsandprelims}

We recall that a \textit{paratopological group} is a pair $(G, \mathcal{T})$ where $G$ is a group and $\mathcal{T}$ is a topology on $G$ such that the mapping $(x, y) \mapsto xy$ of $G\times G$ into $G$ is continuous. If in addition the mapping $x \mapsto x^{-1}$ of $G$ into $G$ is continuous then $(G, \mathcal{T})$ is a topological group.

We call $d: X \times X\to [0, \infty)$ a \textit{quasi-pseudometric} on $X$ if $d(x, y) = 0$ whenever $x = y$ and $d(x, y)\leq d(x, z)+d(z, y)$, for all $x, y, z \in X$. If $d$ is a quasi-pseudometric on a group~$G$ and $d(ax, ay) = d(x, y)$ for all $a, x, y \in G$, then we say that $d$ is \textit{left invariant}; similarly, if $d(xa, ya) = d(x, y)$ for all $a, x, y$, then $d$ is \textit{right invariant}. If $d$ is both left and right invariant, then we say it is \textit{two-sided invariant}. It is easy to check that $d$ is two-sided invariant if and only if
\[
d(x_1 \ldots x_n, y_1 \ldots y_n) \leq d(x_1, y_1) + \cdots + d(x_n, y_n)
\]
for all $x_1, \ldots, x_n, y_1, \ldots, y_n \in G$.

Given a group~$G$ with identity element~$e$, a function $N:G \to [0, \infty)$ is called a \textit{quasi-prenorm} on~$G$ if the following conditions are satisfied:
\begin{enumerate}
\item[(1)]
$N(e)=0$; and
\item[(2)]
$N(gh)\leq N(g)+N(h)$ for all $g, h\in G$.
\end{enumerate}
If $N$ in addition satisfies
\begin{enumerate}
\item[(3)]
$N(h^{-1}gh)=N(g)$ for all $g, h\in G$,
\end{enumerate}
then we say that $N$ is \textit{invariant}.

Let $G$ be a group. If $d$ is a left invariant quasi-pseudometric on~$G$ then the function $N_d:G\to [0, \infty)$ defined by $N_d(x) = d(e, x)$ for all $x \in G$ is a quasi-prenorm on~$G$, and conversely if $N$ is a quasi-prenorm on~$G$ then the function $d_N:G \times G \to [0, \infty)$ defined by $d_N(x, y) = N(x^{-1}y)$ for all $x, y \in G$ is a left invariant quasi-pseudometric on~$G$.
Clearly, the mappings $d \mapsto N_d$ and $N \mapsto d_N$ define a one-to-one correspondence between the family of left invariant quasi-pseudometrics (resp., two-sided invariant quasi-pseudometrics) on~$G$ and the family of quasi-prenorms (resp., invariant quasi-prenorms) on~$G$.

\begin{definition}
\label{Def1}
\begin{rm}
Let $X$ be a subspace of a paratopological group~$G$. Suppose that
\begin{enumerate}
\item[(1)]
the set $X$ generates $G$ algebraically, that is, $\langle X\rangle=G$ and
\item[(2)]
every continuous mapping $f:X\to H$ of $X$ to an arbitrary 
paratopological group $H$ extends to a continuous homomorphism $\hat{f}:G\to H$.
\end{enumerate}
Then $G$ is called the \textit{Markov free paratopological group on~$X$}, and is denoted by $\FP(X)$.

By substituting ``abelian paratopological group'' for each occurrence of ``paratopological group'' above we obtain the definition of the \textit{Markov free abelian paratopological group on~$X$}, which is denoted by $\AP(X)$.
\end{rm}
\end{definition}

\begin{definition}
\label{Def1.5}
\begin{rm}
Let $X$ be a subspace of a paratopological group $G$ and let $e\in X$ be the identity of $G$. Suppose that
\begin{enumerate}
\item[(1)]
$X$ algebraically generates $G$, that is, $\langle X\rangle=G$ and
\item[(2)]
every continuous mapping $f:X\to H$ of $X$ to an arbitrary paratopological group $H$ satisfying $f(e)=e_H$ extends to a continuous homomorphism $\hat{f}:G\to H$.
\end{enumerate}
Then $G$ is called the \textit{Graev free paratopological group on $(X, e)$}, and is denoted by $\FP_G(X, e)$.

By substituting ``abelian paratopological group'' for each occurrence of ``paratopological group'' above we obtain the definition of the \textit {Graev free abelian paratopological group on $(X, e)$}, which is denoted by $\AP_G(X, e)$.
\end{rm}
\end{definition}

\section{Extension of quasi-pseudometrics}

In~\cite{Graev1,Graev2}, Graev developed a method for extending a pseudometric~$d$ from a set~$X$ containing an element~$e$ to a two-sided invariant pseudometric on the abstract free group $F_a(X \setminus \{e\})$ on $X \setminus \{e\}$ (with $e \in X$ identified with the identity element~$e$ of the group), and then employed the method in various applications to Graev free topological groups. A major part of~\cite{RoSaTk} is devoted to the development and application of an analogous process for the extension of a quasi-pseudometric from~$X$ to $F_a(X \setminus \{e\})$. We make substantial use here of ideas and results from~\cite{RoSaTk} relating to this extension process.

Since our applications are to Markov free paratopological groups rather than to Graev free paratopological groups, some changes are required. The changes, however, are fairly minor, and essentially centre around the simple observation that for a topological space~$X$ the groups $\FP(X)$ and $\FP_G(X \oplus \{e\})$ (where `$\oplus$' denotes the topological sum) are topologically isomorphic in a natural way.

We now outline some of the ideas of~\cite{RoSaTk} in a form suitable for our applications. \textit{For most of the remainder of this section, we consider a fixed set~$X$ and a fixed quasi-pseudometric~$d$ on~$X$ which is bounded by~$1$.}

Let $e$ be the identity of the abstract free group $F_a(X)$ on~$X$.
Extend~$d$ from~$X$ to a quasi-pseudometric~$d_e$ on $X\cup\{e\}$ by setting
\[
d_e(x, y)
=
\begin{cases}
      0 & \mbox{if $x=y$}, \\
      d(x, y) & \mbox{if $x, y \in X$}, \\
      1 & \mbox{otherwise}
\end{cases}
\]
for $x, y\in X\cup\{e\}$.
As in~\cite{RoSaTk}, extend~$d_e$ to a quasi-pseudometric~$d^*$ on $\tilde{X} = X\cup\{e\}\cup X^{-1}$ defined by
\[
d^*(x, y)
=
\begin{cases}
      0 & \mbox{if $x=y$}, \\
      d_e(x, y) & \mbox{if $x, y \in X \cup \{e\}$}, \\
      d_e(y^{-1}, x^{-1}) & \mbox{if $x, y \in X^{-1} \cup\{e\}$}, \\
      2 & \mbox{otherwise}
\end{cases}
\]
for $x, y\in \tilde{X}$ (this definition of~$d^*$ is expressed differently from that of~\cite{RoSaTk}, but is easily seen to be equivalent).

\begin{definition}
\label{Def2}
\begin{rm}
Let $H$ be a subset of the set $\N$ of natural numbers such that $|H|=2n$ for some $n\geq 1$. Then a \textit{scheme} \cite{RoSaTk} on~$H$ is a bijection $\varphi: H\to H$ satisfying the following conditions:
\begin{enumerate}
\item[(1)]
if $i\in H$ and $j=\varphi(i)$, then $j\neq i$ and $\varphi(j)=i$; and
\item[(2)]
there are no $i, j\in H$ such that $i<j<\varphi(i)<\varphi(j)$.
\end{enumerate}
We say that $\varphi$ is a \textit{nested scheme} on a set $H=\{i_1, i_2, \cdots, i_{2n}\}\subseteq \N$ where $i_1<i_2<\cdots<i_{2n}$ if $\varphi(i_k)=i_{2n-k+1}$ for all $k=1, 2, \ldots, 2n$.
\end{rm}
\end{definition}

If $\mathcal{X}$ is a word in the alphabet~$\tilde{X}$, then we denote the reduced form of $\mathcal{X}$ by $[\mathcal{X}]$. We denote the length of~$\mathcal{X}$ as a string over~$\tilde{X}$ by~$\ell(\mathcal{X})$.

Let $g\in F_a(X)$ be a reduced word and let $\mathcal{X}$ be a word in the alphabet~$\tilde{X}$ of length $\ell(\mathcal{X}) = 2n$ such that $[\mathcal{X}] = g$. Let $\mathcal{S}_n$ be the family of all schemes $\varphi$ on $\{1, 2, \ldots, 2n\}$. Following~\cite{RoSaTk} we define
\[
\Gamma_d(\mathcal{X}, \varphi)
=
\frac{1}{2}\sum_{i=1}^{2n}d^*(x_i^{-1}, x_{\varphi(i)})
\]
and then we define $N_d: F_a(X) \to [0, \infty)$ by setting $N_d(g) = 0$ if $g = e$ and 
\[
N_d(g)
=
\inf \{\Gamma_d(\mathcal{X}, \varphi):
         [\mathcal{X}] = g, \ \ell(\mathcal{X}) = 2n, \
         \varphi\in \mathcal{S}_n, \ n\in \N\}
\]
for $g \in F_a(X)$ with $g \neq e$. By \cite{RoSaTk}, Claim~3, $N_d$ is an invariant quasi-prenorm on $F_a(X)$. Now let $\hat{d}$ be the two-sided invariant quasi-pseudometric on $F_a(X)$ corresponding to the invariant quasi-prenorm $N_d$ (see section~\ref{defsandprelims}); thus $\hat{d}(g, h) = N_d(g^{-1}h)$ for all $g, h\in F_a(X)$. We refer to~$\hat{d}$ as the \textit{Graev extension} of~$d$ to $F_a(X)$.

\begin{definition}
\label{Def3}
\begin{rm}
If $\mathcal{X}$ is a word in the alphabet $\tilde{X}$, then we say that $\mathcal{X}$ is \textit{almost irreducible} \cite{RoSaTk} if $\mathcal{X}$ does not contain two adjacent symbols $x$ and~$x^{-1}$ for any $x \in \tilde{X}$.
\end{rm}
\end{definition}

\begin{remark}
\label{almost}
\begin{rm}
We note that if $\mathcal{X}$ is an almost irreducible word of length~$2n$, then $\mathcal{X}$ may contain at most $n$ letters equal to~$e$. Also, an almost irreducible word that contains no occurrence of~$e$ is reduced.
\end{rm}
\end{remark}

The following result is essentially Claim~2 of~\cite{RoSaTk}.

\begin{theorem}
\label{claim2}
If $g$ is a reduced word in $F_a(X)$ distinct from $e$, then there exists an almost irreducible word $\mathcal{X}_g = x_1 x_2 \ldots x_{2n}$ of length $2n \geq 2$ in the alphabet~$\tilde{X}$ and a scheme $\varphi_g \in \mathcal{S}_n$ that satisfy the following conditions:
\begin{enumerate}
\item[(1)]
for $i = 1, 2, \ldots, 2n$, either $x_i$ is~$e$ or $x_i$ is a letter in $g$\textup{;}
\item[(2)]
$[\mathcal{X}_g]=g$ and $n \leq \ell (g)$\textup{;} and
\item[(3)]
$N_d(g)=\Gamma_d(\mathcal{X}_g, \varphi_g)$.
\end{enumerate}
\end{theorem}

The next result is probably known, at least in the context of free topological groups, but since we have not found a proof in the literature, we sketch one here. We use the following notation. If $\mathcal{X} = x_1 \ldots x_n$, where $x_1, \ldots, x_n \in \tilde{X}$, then we write $S(\mathcal{X}) = \{x_1, \ldots x_n, x_1^{-1}, \ldots x_n^{-1}\}$.

\begin{theorem}
\label{nested2}
Let $g \in F_a(X)$, let $\mathcal{X}$ be a representation of~$g$ in the alphabet $\tilde{X}$ of length~$2n$ for some $n \geq 1$ and let $\varphi$ be a scheme on the set $\{1, 2, \ldots, 2n\}$. Then there exist a representation~$\mathcal{X}'$ of~$g$ of length~$2m$ for some~$m$ such that $S(\mathcal{X}') = S(\mathcal{X})$ and a nested scheme~$\varphi'$ on $\{1, 2, \ldots, 2m\}$ such that $\Gamma_d(\mathcal{X}, \varphi) = \Gamma_d(\mathcal{X}', \varphi')$.
\end{theorem}

\begin{proof}[Proof outline]
Fix $n \geq 1$ and assume inductively that the desired statement holds for every word in $F_a(X)$, every representation of the word of even length less than~$2n$ and every scheme on the corresponding index set. Consider~$g$, $\mathcal{X}$ and~$\varphi$ as above, and suppose that $\mathcal{X} = x_1 \ldots x_{2n}$.

If $\varphi(1) = 2n$, write $\mathcal{X} = x_1 \mathcal{X}_1 x_{2n}$ and apply the inductive assumption to $\mathcal{X}_1$ and the restriction~$\varphi_1$ of the scheme~$\varphi$ to $\{2, \ldots, 2n-1\}$ (strictly, we should first re-index~$\mathcal{X}_1$ by $\{1, \ldots, 2n-2\}$ and adjust~$\varphi_1$ accordingly). This gives us a word $\mathcal{X}'_1$ and a nested scheme $\varphi'_1$ on a suitable set $\{1, \ldots, 2n'\}$ as in the theorem, and it is clear that we may then construct the desired representation~$\mathcal{X}'$ and nested scheme~$\varphi'$.

Otherwise, there exists~$p$ with $1 \leq p \leq n-1$ such that the restriction of~$\varphi$ to each of the sets $\{1, \ldots, 2p\}$ and $\{2p+1, \ldots, 2n\}$ is a scheme. Write $\mathcal{X} = \mathcal{Y} \mathcal{Z}$, where $\mathcal{Y} = x_1 \ldots x_{2p}$ and $\mathcal{Z} = x_{2p+1} \ldots x_{2n}$,
and apply the inductive assumption to each of $\mathcal{Y}$ and~$\mathcal{Z}$. This gives us respective representations $\mathcal{Y}'$ and~$\mathcal{Z}'$ of lengths $2q$ and $2r$, say, and corresponding nested schemes with the properties in the theorem. Then $\mathcal{Y}' \mathcal{Z}' = x'_1 \ldots x'_{2q} y'_1 \ldots y'_{2r}$ and a scheme~$\psi$ can obviously be constructed from those for $\mathcal{Y}'$ and~$\mathcal{Z}'$ in such a way that the restriction of~$\psi$ to each of $\{1, \ldots, 2q\}$ and $\{2q+1, \ldots, 2q+2r\}$ is nested. Finally, if we define
\[
\mathcal{X}'
 =
x'_1 \ldots x'_{2q}
y'_1 \ldots y'_{2r}
(x'_{2q})^{-1} \ldots (x'_{q+1})^{-1}
x'_{q+1} \ldots x'_{2q}
\]
and let $\varphi'$ be the (unique) nested scheme on $\{1, \ldots, 4q+2r\}$, then it is clear that $\mathcal{X}'$ and~$\varphi'$ have the desired properties. The result follows by induction.
\end{proof}

\begin{theorem}
\label{maximal}
The Graev extension~$\hat{d}$ is the maximal two-sided invariant extension of~$d^*$ from $X\cup \{e\}\cup X^{-1}$ to $F_a(X)$.
\end{theorem}

\begin{proof}
Fix $g, h \in F_a(X)$. Then there exists an almost irreducible representation~$\mathcal{X}$ of $g^{-1}h$, where $\ell(\mathcal{X}) = 2n$ for some~$n$, and a scheme~$\varphi$ on the set $\{1, 2, \ldots, 2n\}$ such that $\hat{d}(g, h) = \hat{d}(e, g^{-1}h) = \Gamma_d(\mathcal{X}, \varphi)$. By Theorem~\ref{nested2} there exists a representation $\mathcal{X}' = x_1 x_2 \ldots x_{2m}$ of $g^{-1}h$ of length~$2m$ for some~$m$ and a nested scheme~$\varphi'$ on the set $\{1, 2, \ldots, 2m\}$ such that $\Gamma_d(\mathcal{X}, \varphi) =\Gamma_d(\mathcal{X}', \varphi')$.

Let $\sigma$ be any two-sided invariant quasi-pseudometric on $F_a(X)$ such that $\sigma|_{\tilde{X}} = d^*$ and write $a = x_1 x_2 \ldots x_m$ and $b  = x_{m+1} x_{m+2} \ldots x_{2m}$. Then
\begin{align*}
\sigma(g, h) 
 &= \sigma(e, g^{-1}h) \\
 &= \sigma(e, ab) \\
 &= \sigma(a^{-1}, b) \\
 &\leq \sum_{i=1}^m \sigma(x_i^{-1}, x_{2m-i+1}) \\
 &= \sum_{i=1}^md^*(x_i^{-1}, x_{2m-i+1}) \\
 &= \Gamma_d (\mathcal{X}', \varphi') \\
 &= \hat{d}(g, h),
\end{align*}
and the result follows.
\end{proof}

For $x \in X$ and $\epsilon > 0$ we denote the ball $\{y: d(x, y) < \epsilon\}$ of radius $\epsilon$ with centre~$x$ by $B_d(x, \epsilon)$.

The next result is Claim~6 of~\cite{RoSaTk}.

\begin{theorem}
\label{RSTtheorem}
The family $\{B_{\hat{d}}(e, \epsilon):\varepsilon>0\}$ is a base at the identity~$e$ for a paratopological group topology~$\mathcal{T}_d$ on the free group $F_a(X)$ and the restriction of $\mathcal{T}_d$ to~$X$ coincides with the topology on~$X$ generated by~$d$.
\end{theorem}

We recall that a real-valued function $f$ on a topological space~$X$ is said to be \textit{upper semi-continuous} if the set $\{x\in X: f(x)<a\}$ is an open set in~$X$ for every $a\in \R$. The \textit{upper topology}~$\tau_u$ for the set~$\R$ has a base of sets of the form $\{x\in \R:x<a\}$ for all $a\in \R$. Clearly, $f$ is upper semi-continuous if and only if $f:X\to (\R, \tau_u)$ is continuous.

If $d$ is a quasi-pseudometric on a space~$X$, then for each $x\in X$ we define $d_x(y)=d(x, y)$ for all $y\in X$. It is easy to see that $d_x$ is upper semi-continuous for all $x\in X$ if and only if the set $B_d(x, \epsilon)$ is open in~$X$ for all $\epsilon>0$.

Let $\mathcal{Q}$ be a family of quasi-pseudometrics on a set $X$ and let
\[
\mathscr{B}
 =
\{B_\rho(x, \epsilon):
    x \in X, \ \rho \in \mathcal{Q} \ \mbox{ and }\ \epsilon > 0\}.
\]
Then we call the topology on~$X$ which has $\mathscr{B}$ as a subbase  \textit{the topology generated by the family $\mathcal{Q}$}.

Every topological space $X$ is generated by a family of quasi-pseudometrics~$\rho$ such that $\rho_x$ is upper semi-continuous for all $x\in X$ (see \cite{Reilly2} and~\cite[page~28]{FlLi}). Specifically, for every open set $U$ in $X$ and for all $x, y\in X$ define $\rho_U$ by
\[
\rho_U(x, y)
=
\begin{cases}
      1 & \mbox{if $x \in U, y \notin U$}, \\
      0 & \mbox{otherwise}.
\end{cases}
\]
Then it is obvious that $\rho_U$ is a quasi-pseudometric on $X$, that $(\rho_U)_x$ is upper semi-continuous for each $x\in X$ and that the family $\mathcal{Q}=\{\rho_U:U$ open in $X\}$ generates the topology of $X$.

Let $X$ be a topological space and let $\mathcal{D}_1$ be the family of all quasi-pseudometrics $d$ on~$X$ which are bounded by~$1$ and are such that $d_x$ is upper semi-continuous for all $x\in X$. Clearly, the family $\mathcal{D}_1$ generates the original topology on~$X$. For every $d \in \mathcal{D}_1$ let $\hat{d}$ be the Graev extension of~$d$ to $F_a(X)$. For each~$d$, Theorem~\ref{RSTtheorem} shows that $\mathcal{T}_d$ is a paratopological group topology on $F_a(X)$ which induces on~$X$ the topology induced by~$d$. It follows that the supremum of all the topologies $\mathcal{T}_d$ for $d\in \mathcal{D}_1$ is a paratopological group topology on $F_a(X)$ which induces the original topology on~$X$. We refer to this topology as the \textit{Graev topology} and denote it by $\mathcal{T}_G$.

Since each topology~$\mathcal{T}_d$ is locally invariant, it follows that the Graev topology is also locally invariant. Now using Proposition~3.1 of~\cite{Ravsky1}, it is easy to see that if $U$ is an open neighbourhood of~$e$ in any locally invariant paratopological group~$G$, then there exists a two-sided invariant quasi-pseudometric~$d$ bounded by~$1$ on~$G$ such that $d_x$ is upper semi-continuous for all $x \in G$ and $B_d(e, 1) \subseteq U$. A straightforward argument using Theorem~\ref{maximal} then yields the following result.

\begin{theorem}
\label{graev_top}
The Graev topology is the finest locally invariant paratopological group topology on $F_a(X)$ which induces the original topology on~$X$. The corresponding topology on the free abelian group $A_a(X)$ is the free topology.
\end{theorem}

\section{Results}

If $Y$ is a subspace of a space~$X$ and $y \in Y$, we write $\cl_Y(y)$ to denote the closure of the singleton~$\{y\}$ in the subspace~$Y$.

\begin{lemma}
\label{rez_lem}
If $G$ is a paratopological group, $A$ is a subset of~$G$ and $a, b \in A$,
then $a \in \cl_A(b)$ if and only if $b^{-1}\in \cl_{A^{-1}}(a^{-1})$.
\end{lemma}

A theorem of Reznichenko (see Theorem~2.4 of Pyrch and Ravsky~\cite{PyRa}) states that if $X$ is an arbitrary topological space and $A$ is a closed subset of~$X$, then $A^{-1}$ is an open subset of the subspace $X^{-1}$ of $\FP(X)$.

Consider the topology $\mathcal{T}_A$ on $X^{-1}$ which has as a base the collection $\{A^{-1}: \mbox{$A$ closed in $X$}\}$; equivalently, $\mathcal{T}_A$ has the collection $\{(\cl_X(x))^{-1} : x \in X\}$ as a base. Clearly, $\mathcal{T}_A$ is closed under arbitrary intersections and each point $x^{-1} \in X^{-1}$ has a smallest open neighbourhood, namely $(\cl_X(x))^{-1}$, and hence $(X^{-1}, \mathcal{T}_A)$ is an Alexandroff space (see~\cite{Arenas}). Moreover, the topology $\mathcal{T}_A$ is the group-theoretical inverse of the so-called Alexandroff dual of the original topology of~$X$ (see Kopperman~\cite{Kopp}).

Let $\mathcal{B}$ be the collection of sets of the form $\{n, n+1, \ldots\} \subseteq \Z$ for all $n\in \Z$. Then $\mathcal{B}$ is a base for a paratopological group topology on the group~$\Z$ of integers under addition, and we denote the corresponding paratopological group by~$\Z^*$.

We answer some obvious questions raised by Reznichenko's result as follows.

\begin{theorem}
\label{rez_thm}
Let $X$ be a topological space. Then $\mathcal{T}_A$, the induced Graev topology~$\mathcal{T}_G|_{X^{-1}}$ and the induced free topology~$\mathcal{T}_F|_{X^{-1}}$ are equal on~$X^{-1}$.
\end{theorem}

\begin{proof}
We show first that $\mathcal{T}_A \subseteq \mathcal{T}_G|_{X^{-1}}$. Let $A$ be a closed subset of~$X$ and, following~\cite{PyRa}, define $f:X \to \Z^*$ by mapping all elements of~$A$ to~$0$ and all other elements of~$X$ to~$1$. Clearly, $f$ is continuous. Extend~$f$ to a homomorphism $\hat{f}:F_a(X) \to \Z^*$. Since $\Z^*$ is an abelian and hence locally invariant paratopological group, Theorem~\ref{graev_top} implies that $\hat{f}$ is continuous with respect to~$\mathcal{T}_G$. Therefore, $A^{-1} = \hat{f}^{-1}(\{0, 1, \ldots\}) \cap X^{-1}$ is open in~$X^{-1}$ with the topology $\mathcal{T}_G|_{X^{-1}}$, and it follows that $\mathcal{T}_A \subseteq \mathcal{T}_G|_{X^{-1}}$.

Clearly, we have $\mathcal{T}_G|_{X^{-1}}\subseteq \mathcal{T}_F|_{X^{-1}}$.

To show that $\mathcal{T}_F|_{X^{-1}}\subseteq \mathcal{T}_A$, consider any fixed $x\in X$. For any $y^{-1}\in (\cl_X(x))^{-1}$ we have $y\in \cl_X(x)$, and Lemma~\ref{rez_lem} implies that $x^{-1}\in \cl_{X^{-1}}(y^{-1})$. If $U$ is a neighbourhood of $x^{-1}$ in $\mathcal{T}_F|_{X^{-1}}$, we therefore have $y^{-1} \in U$, and it follows that $(\cl_X(x))^{-1}\subseteq U$. Thus $\mathcal{T}_F|_{X^{-1}}\subseteq \mathcal{T}_A$, and the proof is complete.
\end{proof}

The following result was noted, in the case of the free topology, in~\cite{PyRa}.

\begin{corollary}
\label{Xinverse}
If $X$ is a $T_1$ space then the Graev topology and the free topology on $X^{-1}$ are discrete.
\end{corollary}

Clearly if $X$ is a $T_1$ space and $x_1, x_2, \ldots, x_n$ are distinct points in $X$, then there exist open sets $U_1, U_2, \ldots, U_n$ in $X$ containing $x_1, x_2, \ldots, x_n$, respectively, such that  $x_i\notin U_j$ whenever $i\neq j$, for $i, j=1, 2, \ldots, n$.

\begin{lemma}
\label{qpsm}
Let $X$ be a $T_1$ space, let $x_1, x_2, \ldots, x_n$ be distinct points of $X$ and suppose that open sets $U_1, U_2, \ldots, U_n$ are chosen as above. Then for $i, j=1, 2, \ldots, n$ with $i\neq j$ the function $d_{i,j}$ defined by setting
\[
d_{i,j}(x, y)
=
\begin{cases}
      1 & \mbox{if $(x\neq x_j \mbox{ and } y = x_j)$} \\
        & \mbox{or $(x \in U_i \mbox{ and } y \notin U_i)$},\\
      0 & \mbox{otherwise}
\end{cases}
\]
for all $x, y\in X$ is a quasi-pseudometric, and $(d_{i,j})_x$ is upper semi-continuous for all $x\in X$.
\end{lemma}

\begin{proof}
It is straightforward to check that $d_{i,j}$ may be represented equivalently by the formula
\[
d_{i,j}(x, y)
=
\begin{cases}
      0 & \mbox{if $x, y \in U_i$} \\
        & \mbox{or $x = x_j$} \\
        & \mbox{or $(x \notin U_i, x\neq x_j \mbox{ and } y \neq x_j)$}, \\
      1 & \mbox{otherwise}
\end{cases}
\]
for all $x, y \in X$; we also observe that the three disjuncts in the first part of this alternative expression are mutually exclusive. Fix $i$ and~$j$ with $i \neq j$. We show first that $d_{i,j}$ is a quasi-pseudometric on~$X$. Clearly, $d_{i,j}(x, y)\geq 0$ for all $x, y \in X$. To show that $d_{i,j}(x, y) \leq d_{i,j}(x, z) + d_{i,j}(z, y)$ for $x, y, z \in X$, it obviously suffices to consider the case when $d_{i,j}(x, y) = 1$. There are two sub-cases. First, suppose that $x \neq x_j$ and $y = x_j$. If $d_{i,j}(x, z) = 0$ then either $x, z \in U_i$, which gives $d_{i,j}(z, y) = 1$, or $x  \notin U_i$ and $z \neq x_j$, again giving $d_{i,j}(z, y) = 1$. Second, suppose that $x \in U_i$ and $y \notin U_i$. If $d_{i,j}(x, z) = 0$ it follows that $z \in U_i$, which gives $d_{i,j}(z, y) = 1$. Therefore the triangle inequality holds.

To show that $(d_{i,j})_x$ is upper semi-continuous, consider $x \in X$ and $\epsilon > 0$. If $x \in U_i$, then $B_{d_{i,j}}(x, \epsilon) = \{y: d_{i,j}(x, y) < \epsilon\} = U_i$ when $\epsilon \leq 1$ and $B_{d_{i,j}}(x, \epsilon) = X$ when $\epsilon > 1$. If $x \notin U_i$ and $x \neq x_j$, then $B_{d_{i,j}}(x, \epsilon) = X \setminus \{x_j\}$ when $\epsilon \leq 1$ and $B_{d_{i,j}}(x, \epsilon) = X$ when $\epsilon > 1$. If $x = x_j$, then $B_{d_{i,j}}(x, \epsilon) = X$. Therefore, $d_{i,j}$ is upper semi-continuous.
\end{proof}

\begin{lemma}
\label{dxixj}
With hypotheses and notation as above, we have the following for $i \neq j$.
\begin{enumerate}
\item
$d_{i,j}(x_i, x) = 0$ if and only if $x \in U_i$.
\item
$d_{i,j}(x, x_j) = 0$ if and only if $x = x_j$.
\item
$d_{i,j}(x_i, x_j) = 1$.
\end{enumerate}
\end{lemma}

Now we state and prove our main theorem.

\begin{theorem}
\label{joiner}
Let $X$ be a $T_1$ space and let $w=x_1^{\epsilon_1}x_2^{\epsilon_2}\ldots x_n^{\epsilon_n}$ be a reduced word in $\FP_n(X)$, where $x_i\in X$ and $\epsilon_i=\pm 1$ for $i=1, 2, \ldots, n$ and if $x_i=x_{i+1}$ for some $i=1, 2, \ldots, n-1$ then $\epsilon_i=\epsilon_{i+1}$. Let $\mathscr{B}$ denote the collection of all sets of the form $U_1^{\epsilon_1} U_2^{\epsilon_2} \ldots U_n^{\epsilon_n}$, where for $i=1, 2, \ldots, n$ the set~$U_i$ is a neighbourhood of $x_i$ in $X$ when $\epsilon_i=1$ and $U_i = \{x_i\}$ when $\epsilon_i = -1$. Then $\mathscr{B}$ is a base for the neighbourhood system at~$w$ in the subspace $\FP_n(X)$ of $\FP(X)$.
\end{theorem}

\begin{proof}
(1) We show that every neighbourhood of~$w$ in $\FP_n(X)$ contains an element of the collection~$\mathscr{B}$. Let $W$ be a such neighbourhood, so that $W = V \cap \FP_n(X)$ for some neighbourhood~$V$ of~$w$ in $\FP(X)$. Since $\FP(X)$ is a paratopological group, there exist in $\FP(X)$ neighbourhoods $V_1, V_2, \ldots, V_n$ of $x_1^{\epsilon_1}, x_2^{\epsilon_2}, \ldots, x_n^{\epsilon_n}$, respectively, such that $w \in V_1 V_2 \ldots V_n \subseteq V$. When $\epsilon_i = 1$ let $U_i = V_i \cap X$ and when $\epsilon_i = -1$ let $U_i = \{x_i\}$. Then $U_1^{\epsilon_1} U_2^{\epsilon_2} \ldots U_n^{\epsilon_n} \subseteq V_1 V_2 \ldots V_n$, and so, setting $B = U_1^{\epsilon_1} U_2^{\epsilon_2} \ldots U_n^{\epsilon_n}$, we have $B \in \mathscr{B}$ and $w \in B \subseteq W$, as required.

\medskip
(2) We show that every element of~$\mathscr{B}$ is a neighbourhood of~$w$ in $\FP_n(X)$. Thus, for $i = 1, 2, \ldots, n$ we suppose that $U_i$ is a fixed neighbourhood of $x_i$ if $\epsilon_i=1$ and that $U_i = \{x_i\}$ when $\epsilon_i = -1$, and we consider $B = U_1^{\epsilon_1} U_2^{\epsilon_2}\ldots U_n^{\epsilon_n} \in \mathscr{B}$.

Choose indices $i_1, i_2, \ldots, i_{n_1}$ for some $n_1\leq n$ such that $x_{i_1}, x_{i_2}, \ldots, x_{i_{n_1}}$ are the distinct letters among $x_1, x_2, \ldots, x_n$, and write $A=\{1, 2, \ldots, n_1\}$. For each $j \in A$, define
\[
I_j = \{i: 1 \leq i \leq n, \ x_i = x_{i_j} \mbox{ and } \epsilon_i = 1\}.
\]
Now pick open neighbourhoods $V_1, V_2, \ldots, V_{n_1}$ of $x_{i_1}, x_{i_2}, \ldots, x_{i_{n_1}}$ in~$X$, respectively, such that
\begin{enumerate}
\item[(i)]
for all $j \in A$, we have $V_j \subseteq U_i$ for all $i \in I_j$, and
\item[(ii)]
for all $j, k \in A$ with $j \neq k$, we have $x_{i_k} \notin V_j$.
\end{enumerate}
For each $j, k \in A$ with $j \neq k$, define
\[
d_{j,k}(x, y)
=
\begin{cases}
      1 & \mbox{if $(x \neq x_{i_k} \mbox{ and } y = x_{i_k})$} \\
        & \mbox{or $(x \in V_j \mbox{ and } y \notin V_j)$},\\
      0 & \mbox{otherwise}.
\end{cases}
\]
By Lemma~\ref{qpsm}, each~$d_{j,k}$ is a quasi-pseudometric on~$X$ such that $(d_{j,k})_x$ is upper semi-continuous for each $x \in X$. Hence, if we define
\[
d(x, y) = \max\{d_{j,k}(x, y): j, k \in A \mbox{ and } j \neq k\},
\]
then $d$ is also a quasi-pseudometric on~$X$ such that $d_x$ is upper semi-continuous for each $x \in X$. Let $\hat{d}$ be the Graev extension of~$d$ to $F_a(X)$. We will show that $B$ is a neighbourhood of~$w$ in $\FP_n(X)$ by showing that $B_{\hat{d}}(e, 1) w \cap \FP_n(X)\subseteq B$.

Let
\[
h = y_1^{\delta_1} y_2^{\delta_2} \ldots y_{p}^{\delta_{p}}
  \in B_{\hat{d}}(e, 1) w \cap \FP_n(X)
\]
be a reduced word of length $\ell(h) = p \leq n$. Then
\[
hw^{-1}
 =
y_1^{\delta_1} y_2^{\delta_2} \ldots y_{p}^{\delta_{p}}
x_n^{-\epsilon_n} x_{n-1}^{-\epsilon_{n-1}} \ldots x_1^{-\epsilon_1}
\in B_{\hat{d}}(e, 1).
\]
Although $h$ and $w$ are reduced, cancellation may occur in the product $hw^{-1}$. Assume that the number of cancelling pairs in $hw^{-1}$ is $\alpha$, where $0 \leq \alpha \leq p$, so that $y_{p-\beta+1} = x_{n-\beta+1}$ and $\delta_{p-\beta+1} = \epsilon_{n-\beta+1}$ for $\beta = 1, 2, \ldots, \alpha$. Write $g = hw^{-1}$, so that in reduced form we have
\[
g
 =
y_1^{\delta_1} y_2^{\delta_2}\ldots y_l^{\delta_l}
x_m^{-\epsilon_m} x_{m-1}^{-\epsilon_{m-1}} \ldots x_1^{-\epsilon_1},
\]
where $l = p-\alpha$ and $m = n-\alpha$. Since $\hat{d}(e, g)<1$, we have $\hat{d}(e, g)=N_d(g)=0$.

If $g = e$ then $h = w \in B$ and there is nothing to prove, so let us assume that $g \neq e$. Then by Theorem~\ref{claim2}, there exist an almost irreducible word $\mathcal{Z}_g = z_1 z_2 \ldots z_{2m_1}$ for some $m_1 \geq 1$ and a scheme $\varphi_g$ on the set $H_1 = \{1, 2, \ldots, 2m_1\}$ such that
\begin{enumerate}
\item[(i)]
each $z_i$ is either~$e$ or a letter in~$g$,
\item[(ii)]
$[\mathcal{Z}_g] = g$ and $\ell(\mathcal{Z}_g) = 2m_1 \leq 2\ell(g) = 2(l+m)$, and
\item[(iii)]
$N_d(g) = \Gamma_d(\mathcal{Z}_g, \varphi_g) = 0$.
\end{enumerate}
From~(iii), we have
\[
\Gamma_d(\mathcal{Z}_g, \varphi_g)
 = \frac{1}{2}\sum_{i=1}^{2m_1} d^*(z_i^{-1}, z_{\varphi_g(i)})
 = 0,
\]
and so
\begin{equation}
\label{eqn10}
d^*(z_i^{-1}, z_{\varphi_g(i)}) = 0 \mbox{ for all $i\in H_1$}.
\end{equation}
Now if $z_i=e$ for any $i\in H_1$ then also $z_{\varphi_g(i)} = e$, because if $z_{\varphi_g(i)} \neq e$ then $d^*(z_i^{-1}, z_{\varphi_g(i)}) = 1$ by definition of~$d^*$, which is impossible by~(\ref{eqn10}).

If all occurrences of~$e$ are removed from~$\mathcal{Z}_g$, then by Remark~\ref{almost} the resulting word~$\mathcal{Z}'_g$ is reduced, so that $\mathcal{Z}'_g = g$ and $\mathcal{Z}'_g$ in particular has length $l+m$. Let us write
\[
\mathcal{Z}'_g = z'_1 z'_2 \ldots z'_{l+m}.
\] 
Moreover, since $d^*(e, e) = 0$, we may use the scheme~$\varphi_g$ on $H_1 = \{1, 2, \ldots, 2m_1\}$ to define a scheme~$\varphi'_g$ on $H_2 = \{1, 2, \ldots, l+m\}$ with the property that $N_d(g) = \Gamma_d(\mathcal{Z}'_g, \varphi'_g) = 0$. Formally, suppose that when the indices among the elements of~$H_1$ corresponding to occurrences of~$e$ in $\mathcal{Z}_g$ are removed, the indices remaining form the set $J = \{j_1, j_2, \ldots, j_{l+m}\}$, where $j_1 < j_2 < \cdots < j_{l+m}$.
Now let $H_2 = \{1, 2, \ldots, l+m\}$ and let $f: J \to H_2$ be the bijection given by $f(j_k) = k$ for $k = 1, 2, \ldots, l+m$. Then it is easy to check that the map $\varphi'_g: H_2 \to H_2$ defined by $\varphi'_g(k) = f(\varphi_g(f^{-1}(k)))$ for $k \in H_2$ is a scheme on~$H_2$ and has the properties claimed.

Let us now for convenience suppress the prime superscripts used above, so that we have
\[
\mathcal{Z}_g
 =
z_1 z_2 \ldots z_{l+m}
 =
y_1^{\delta_1} y_2^{\delta_2}\ldots y_l^{\delta_l}
x_m^{-\epsilon_m} x_{m-1}^{-\epsilon_{m-1}} \ldots x_1^{-\epsilon_1}
 =
g
\]
and $\varphi_g$ is a scheme on~$H_2$ such that $\Gamma_d(\mathcal{Z}_g, \varphi_g) = 0$. From the last equation, we have
\begin{equation}
\label{eqn20}
d^*(z_i^{-1}, z_{\varphi_g(i)}) = 0 \mbox{ for all $i \in H_2$}.
\end{equation}

We claim now that $l = m$ and hence that $p = n$. Assume that $l < m$. Then there exist $q \geq 1$ and distinct $k_1, \ldots, k_q, l_1, \ldots, l_q \in H_2$ such that $l_r = \varphi_g(k_r)$ and $l+1 \leq k_r, l_r \leq l+m$ for $r = 1, 2, \ldots, q$. For any $r = 1, 2, \ldots, q$, set $s \equiv s(r) = l+m+1-k_r$ and $t \equiv t(r) = l+m+1-l_r$, so that 
\[
z_{k_r} = x_s^{-\epsilon_s}
 \quad\mbox{and}\quad
z_{l_r} = x_t^{-\epsilon_t}.
\]
This gives
\begin{equation}
\label{eqn30}
d^*(z_{k_r}^{-1}, z_{\varphi_g(k_r)})
 = d^*(z_{k_r}^{-1}, z_{l_r})
 = d^*(x_s^{\epsilon_s}, x_t^{-\epsilon_t}).
\end{equation}
If $\epsilon_s = \epsilon_t$, then either $\epsilon_s = \epsilon_t = 1$, and we have
\[
d^*(x_s^{\epsilon_s}, x_t^{-\epsilon_t})
 = d^*(x_s, x_t^{-1})
 > 0,
\]
or $\epsilon_s = \epsilon_t = -1$, and we have
\[
d^*(x_s^{\epsilon_s}, x_t^{-\epsilon_t})
 = d^*(x_s^{-1}, x_t)
 > 0,
\]
and in both cases we conclude from~(\ref{eqn30}) that $d^*(z_{k_r}^{-1}, z_{\varphi_g(k_r)}) > 0$, which contradicts~(\ref{eqn20}). Therefore, for $r = 1, 2, \ldots, q$, we have $\epsilon_s = -\epsilon_t$.

For any~$r$ such that $\epsilon_s \equiv \epsilon_{s(r)}= 1$ and $\epsilon_t \equiv \epsilon_{t(r)} = -1$, we find from (\ref{eqn20}) and~(\ref{eqn30}) that
\[
d^*(x_s^{\epsilon_s}, x_t^{-\epsilon_t})
 = d^*(x_s, x_t)
 = 0
\]
and hence that $d_{j,k}(x_s, x_t) = 0$ for all $j, k \in A$ with $j \neq k$,
while if $\epsilon_s = -1$ and $\epsilon_t = 1$, we find that
\[
d^*(x_s^{\epsilon_s}, x_t^{-\epsilon_t})
 = d^*(x_s^{-1}, x_t^{-1})
 = d^*(x_t, x_s)
 = 0
\]
and hence that $d_{j,k}(x_t, x_s) = 0$ for all $j, k \in A$ with $j \neq k$.
Therefore, in either case, Lemma~\ref{dxixj} part~(3) shows that $x_s = x_t$.

Pick~$r$ so that $|s - t| \equiv |s(r) - t(r)| = |k_r - l_r|$ is minimal. Now $|s - t|$ cannot equal~$1$, since the fact that $x_s = x_t$ and $\epsilon_s = -\epsilon_t$ would then contradict the hypothesis that the word~$w$ is reduced. Therefore, by the definition of a scheme, there exists~$r'$ such that $k_r < k_{r'}, l_{r'} < l_r$ or $l_r < k_{r'}, l_{r'} < k_r$, and this contradicts the minimality of $|s - t|$.

This contradiction implies that $l = m$, from which it follows immediately that $p = n$. Furthermore, the argument above shows that if $m+1 \leq i \leq 2m$ then $1 \leq \varphi_g(i) \leq m$ and if $1 \leq i \leq m$ then $m+1 \leq \varphi_g(i) \leq 2m$, for all $i \in H_2$. It follows that $\varphi_g(1) = 2m$, because if $\varphi_g(1) < 2m$ then the fact that $\varphi_g$ is a scheme on $H_2 = \{1, 2, \ldots, 2m\}$ would imply that there exist $i, j \in H_2$ with $\varphi_g(1) < i, j \leq 2m$ such that $\varphi_g(i) = j$ and $\varphi_g(j) = i$, contradicting what we have just shown. Continuing similarly, we find that $\varphi_g(i) = 2m-i+1$ for all $i \in H_2$, that is, that $\varphi_g$ is a nested scheme on~$H_2$.

Therefore,
\[
\Gamma_d(\mathcal{Z}_g, \varphi_g)
 = \sum_{i=1}^{m}d^*(y_i^{-\delta_i}, x_i^{-\epsilon_i})
 = 0,
\]
and so $d^*(y_i^{-\delta_i}, x_i^{-\epsilon_i}) = 0$ for $i = 1, 2, \ldots, m$. It follows that $\delta_i = \epsilon_i$ for all $i = 1, 2, \ldots, m$.
If $\epsilon_i = 1$ for any~$i$, then $d^*(y_i^{-\delta_i}, x_i^{-\epsilon_i}) = d^*(x_i, y_i) = 0$ and hence $d_{j,k}(x_i, y_i) = 0$ for all $j, k \in A$ with $j \neq k$. But there exists $j_0 \in A$ such that $x_{i_{j_0}} = x_i$, so it follows by Lemma~\ref{dxixj} part~(1) that $y_i \in V_{j_0} \subseteq U_i$.
If $\epsilon_i = -1$ for any~$i$, then $d^*(y_i^{-\delta_i}, x_i^{-\epsilon_i}) = d^*(y_i, x_i) = 0$ and hence $d_{j,k}(y_i, x_i) = 0$ for all $j, k \in A$ with $j \neq k$. But there exists $k_0 \in A$ such that $x_{i_{k_0}} = x_i$, so it follows by Lemma~\ref{dxixj} part~(2) that $y_i = x_{i_{k_0}}= x_i$.

Finally,
\begin{align*}
h
 &= y_1^{\delta_1} y_2^{\delta_2} \ldots y_p^{\delta_p} \\
 &= y_1^{\delta_1} y_2^{\delta_2} \ldots y_n^{\delta_n} \\
 &= y_1^{\epsilon_1} y_2^{\epsilon_2} \ldots y_m^{\epsilon_m} y_{m+1}^{\delta_{m+1}} \ldots y_n^{\delta_n}\\
 &\in
U_1^{\epsilon_1} U_2^{\epsilon_2} \ldots U_m^{\epsilon_m} y_{m+1}^{\delta_{m+1}} \ldots y_n^{\delta_n} \\
 &=
U_1^{\epsilon_1} U_2^{\epsilon_2} \ldots U_m^{\epsilon_m}x_{m+1}^{\epsilon_{m+1}} \ldots x_n^{\epsilon_n} \\
 &\subseteq
U_1^{\epsilon_1} U_2^{\epsilon_2} \ldots U_n^{\epsilon_n}.
\end{align*}
Therefore, $B_{\hat{d}}(e, 1) w \cap \FP_n(X)\subseteq U_1^{\epsilon_1} U_2^{\epsilon_2} \ldots U_n^{\epsilon_n} = B$, and so $B$ is a neighbourhood of~$w$ in $\FP_n(X)$, as required.
\end{proof}

\begin{remark}
\label{joinerremark1}
\begin{rm}
Part (1) of the proof of Theorem~\ref{joiner} remains valid for any paratopological group topology on $F_a(X)$ that induces the original topology on~$X$, so it follows that $\mathscr{B}$ is a base for the neighbourhood system at~$w$ in the subspace $\FP_n(X)$ of $F_a(X)$ when the latter is equipped with the Graev topology~$\mathcal{T}_G$.
\end{rm}
\end{remark}

\begin{remark}
\label{joinerremark2}
\begin{rm}
It is clear from the proof of Theorem~\ref{joiner} that for each $B \in \mathscr{B}$ there exists $B' \in \mathscr{B}$ such that $B' \subseteq B$ and every element of~$B'$ is of reduced length exactly~$n$.
\end{rm}
\end{remark}

The analogue of Theorem~\ref{joiner} for the free abelian paratopological group takes the following form; the proof is similar to the proof above, and is omitted.

\begin{theorem}
Let $X$ be a $T_1$ space and let $w = \epsilon_1 x_1 + \epsilon_2 x_2 + \cdots + \epsilon_n x_n$ be a reduced word in $\AP_n(X)$, where $x_i\in X$ and $\epsilon_i=\pm 1$ for all $i=1, 2, \ldots, n$ and if $x_i=x_j$ for some $i, j=1, 2, \ldots, n$ then $\epsilon_i=\epsilon_j$. Then the collection $\mathscr{B}$ of all sets of the form $\epsilon_1U_1+\epsilon_2U_2+\cdots +\epsilon_nU_n$, where for all $i=1, 2, \ldots, n$ the set $U_i$ is a neighbourhood of $x_i$ in $X$ when $\epsilon_i=1$ and $U_i=\{x_i\}$ when $\epsilon_i=-1$, is a base for the neighbourhood system at $w$ in $\AP_n(X)$.
\end{theorem}

The following result was proved in Proposition~3.9 in \cite{PyRa} under the stronger hypothesis that $X$ is Tychonoff. Given Theorem~\ref{joiner}, the proof is essentially identical to the proof of the corresponding result for free topological groups given in~\cite{HaMoTh}.

\begin{theorem}
\label{Xpowers}
Let $X$ be a $T_1$ space. Then the free paratopological group $\FP(X)$ contains as a closed subspace a homeomorphic copy of the product space $X^n$ for each $n \geq 1$.
\end{theorem}

A result similar to the following was given in~\cite{PyRa} for the case of the free abelian paratopological group $\AP(X)$. 

\begin{theorem}
\label{equiv_conds}
The following conditions are equivalent for a topological space~$X$.
\begin{enumerate}
\item[(1)]
The space $X$ is $T_1$.
\item[(2)]
The space $\FP(X)$ is $T_1$.
\item[(3)]
The subspace $X$ of $\FP(X)$ is closed.
\item[(4)]
The subspace $X^{-1}$ of $\FP(X)$ is discrete.
\item[(5)]
The subspace $X^{-1}$ of $\FP(X)$ is $T_1$.
\item[(6)]
The subspace $X^{-1}$ of $\FP(X)$ is closed.
\item[(7)]
The subspace $\FP_n(X)$ of $\FP(X)$ is closed for all $n \in \N$.
\item[(8)]
The subspace $\FP_n(X)$ of $\FP(X)$ is closed for some $n \in \N$.
\end{enumerate}
\end{theorem}

\begin{proof}
A convenient scheme of proof is to show that
(1) $\Rightarrow$ (3)~$\Rightarrow$ (2)~$\Rightarrow$~(1),
(1) $\Rightarrow$ (4)~$\Rightarrow$ (5)~$\Rightarrow$~(1),
(1) $\Rightarrow$ (6)~$\Rightarrow$~(2) and
(1) $\Rightarrow$ (7)~$\Rightarrow$ (8)~$\Rightarrow$~(2).
However, the only implications here that are not either trivial or given by rewriting arguments from the corresponding proof in~\cite{PyRa} in non-abelian notation are those for (1)~$\Rightarrow$~(3), (1)~$\Rightarrow$~(6) and (1)~$\Rightarrow$~(7), so we prove only these (the argument from~\cite{PyRa} can also be adapted to show that (1)~$\Rightarrow$~(6), but we give a simpler proof).

First, for each $n \in \N$, let $Z_n$ be the subset of $\FP(X)$ consisting of the words of exponent sum~$n$. Then $Z_n$ is open, since $Z_n = \hat{f}^{-1}(\{n\})$ where $\hat{f}: \FP(X) \to \Z$ is the continuous homomorphism extending the continuous function $f: X \to \Z$ defined by $f(x) = 1$ for all $x \in X$.

Now assume that $X$ is $T_1$.

(1) $\Rightarrow$ (3):
To show that $X$ is closed in $\FP(X)$, let $w$ be a reduced word in $\FP(X)$ such that $w \notin X$.  If $w \in \FP_1(X)$, then either $w \in X^{-1} \subseteq Z_{-1}$ or $w = e \in Z_0$, and $Z_{-1}$ and~$Z_0$ are open and disjoint from $X \subseteq Z_1$. If $w \notin \FP_1(X)$, let $n \geq 1$ be the smallest natural number such that $w \notin \FP_n(X)$. Then $w \in \FP_{n+1}(X) \setminus \FP_n(X)$ and $w$ has length exactly $n+1$. By Theorem~\ref{joiner} and Remark~\ref{joinerremark2} there exists a neighbourhood~$U$ of~$w$ in $\FP_{n+1}(X)$ such that $U \subseteq \FP_{n+1}(X) \setminus \FP_n(X)$. Hence there exists a neighbourhood~$V$ of~$w$ in $\FP(X)$ such that $U = V \cap \FP_{n+1}(X)$ and $V \cap \FP_n(X)= \emptyset$. In particular, $V \cap X = \emptyset$. Therefore, $X$ is closed in $\FP(X)$.

(1) $\Rightarrow$ (7):
Fix $n \in \N$. Let $w \notin \FP_n(X)$ and suppose that $w$ has reduced length~$k > n$. By Theorem~\ref{joiner} and Remark~\ref{joinerremark2} there exists a neighbourhood~$U$ of~$w$ in $\FP_k(X)$ such that $U \subseteq \FP_k(X) \setminus \FP_{k-1}(X)  \subseteq \FP_k(X) \setminus \FP_n(X)$. Hence there exists a neighbourhood~$V$ of~$w$ in $\FP(X)$ such that $U = V \cap \FP_k(X)$ and $V \cap \FP_n(X) = \emptyset$. Therefore, $\FP_n(X)$ is closed in $\FP(X)$.

(1) $\Rightarrow$ (6):
Since (7) holds, $\FP_1(X)$ is closed in $\FP(X)$. But $X^{-1} = Z_{-1} \cap \FP_1(X)$, so $X^{-1}$ is closed in $\FP_1(X)$. Therefore, $X^{-1}$ is closed in $\FP(X)$.
\end{proof}

\def\cprime{$'$} \def\cprime{$'$} \def\cprime{$'$} \def\cprime{$'$}
  \def\cprime{$'$}
\providecommand{\bysame}{\leavevmode\hbox to3em{\hrulefill}\thinspace}


\begin{thebibliography}{10}

\bibitem{Arenas}
F.~G. Arenas, \emph{Alexandroff spaces}, Acta Math. Univ. Comenian. (N.S.)
  \textbf{68} (1999), 17--25.

\bibitem{Arh1}
A.~V. Arhangel{\cprime}ski{\u\i}, \emph{Mappings connected with topological
  groups}, Dokl. Akad. Nauk SSSR \textbf{181} (1968), 1303--1306.

\bibitem{Arh2}
\bysame, \emph{On the relations between invariants of topological groups and
  their subspaces}, Uspekhi Mat. Nauk \textbf{35} (1980), 3--22, International
  Topology Conference (Moscow State Univ., Moscow, 1979).

\bibitem{FlLi}
Peter Fletcher and William~F. Lindgren, \emph{Quasi-uniform spaces}, Lecture
  Notes in Pure and Applied Mathematics, vol.~77, Marcel Dekker Inc., New York,
  1982.

\bibitem{Graev1}
M.~I. Graev, \emph{Free topological groups}, Izvestiya Akad. Nauk SSSR. Ser.
  Mat. \textbf{12} (1948), 279--324.

\bibitem{Graev2}
\bysame, \emph{Free topological groups}, Amer. Math. Soc. Translation
  \textbf{1951} (1951), 61.

\bibitem{HaMoTh}
J.~P.~L. Hardy, Sidney~A. Morris, and H.~B. Thompson, \emph{Applications of the
  {S}tone-\v {C}ech compactification to free topological groups}, Proc. Amer.
  Math. Soc. \textbf{55} (1976), 160--164.

\bibitem{Joiner}
Charles Joiner, \emph{Free topological groups and dimension}, Trans. Amer.
  Math. Soc. \textbf{220} (1976), 401--418.

\bibitem{Kopp}
Ralph Kopperman, \emph{Asymmetry and duality in topology}, Topology Appl.
  \textbf{66} (1995), 1--39.

\bibitem{MaRo}
Josefa Marin and Salvador Romaguera, \emph{A bitopological view of
  quasi-topological groups}, Indian J. Pure Appl. Math. \textbf{27} (1996),
  393--405.

\bibitem{Markov1}
A.~Markov, \emph{On free topological groups}, C. R. (Doklady) Acad. Sci. URSS
  (N. S.) \textbf{31} (1941), 299--301.

\bibitem{Markov2}
A.~A. Markov, \emph{On free topological groups}, Bull. Acad. Sci. URSS. S\'er.
  Math. [Izvestia Akad. Nauk SSSR] \textbf{9} (1945), 3--64.

\bibitem{Markov3}
\bysame, \emph{Three papers on topological groups: {I}. {O}n the existence of
  periodic connected topological groups. {II}. {O}n free topological groups.
  {III}. {O}n unconditionally closed sets}, Amer. Math. Soc. Translation
  \textbf{1950} (1950), 120.

\bibitem{PyRa}
N.~M. Pyrch and O.~V. Ravsky, \emph{On free paratopological groups}, Mat. Stud.
  \textbf{25} (2006), 115--125.

\bibitem{Ravsky1}
O.~V. Ravsky, \emph{Paratopological groups. {I}}, Mat. Stud. \textbf{16}
  (2001), 37--48.

\bibitem{Ravsky2}
\bysame, \emph{Paratopological groups. {II}}, Mat. Stud. \textbf{17} (2002),
  93--101.

\bibitem{Reilly2}
Ivan~L. Reilly, \emph{On generating quasi uniformities}, Math. Ann.
  \textbf{189} (1970), 317--318.

\bibitem{RoSaTk}
S.~Romaguera, M.~Sanchis, and M.~Tkachenko, \emph{Free paratopological groups},
  Proceedings of the 17th {S}ummer {C}onference on {T}opology and its
  {A}pplications, vol.~27, 2003, pp.~613--640.

\end{thebibliography}
\end{document}